\documentclass[11pt,a4paper]{article}

\usepackage[english]{babel}
\usepackage{graphicx} 
\usepackage[dvipsnames]{xcolor}
\usepackage{tikz}
\usetikzlibrary{arrows.meta,backgrounds,calc,fit,positioning,shapes.geometric}
\usepackage{amssymb,amsmath,amscd,amsthm,amsfonts,mathrsfs}
\usepackage{algorithm}
\usepackage{enumitem}
\usepackage{authblk}
\usepackage{hyperref,cite}
\hypersetup{
    colorlinks = true,
    linkcolor = Magenta,
    citecolor = ForestGreen,
    urlcolor = blue,
    linkbordercolor = {blue}
}

\usepackage{comment,todonotes}
\usepackage{lineno}

\usepackage{subcaption}
\usepackage{cleveref}

\tikzset{
  edge/.style={line width=0.55pt},
  marked vertex/.style={line width=0.85pt,line cap=round},
  vertex label/.style={font=\small,inner sep=1pt},
  case label/.style={anchor=west,inner sep=0pt}
}

\newcommand{\Vertex}[2]{%
  \if x#2\relax
    \draw[marked vertex] (#1)
      ++(-1.8pt,-1.8pt) -- ++(3.6pt,3.6pt)
      (#1) ++(-1.8pt,1.8pt) -- ++(3.6pt,-3.6pt);
  \else
    \fill (#1) circle[radius=1.05pt];
  \fi
}

\newcommand{\Comb}[1]{%
  \foreach \lowerMark/\upperMark [count=\i from 0] in {#1}{%
    \coordinate (lower-\i) at (\i,0);
    \coordinate (upper-\i) at (\i,1);
    \draw[edge] (lower-\i) -- (upper-\i);
    \ifnum\i>0
      \pgfmathtruncatemacro{\previous}{\i-1}
      \draw[edge] (lower-\previous) -- (lower-\i);
    \fi
  }%
  \foreach \lowerMark/\upperMark [count=\i from 0] in {#1}{%
    \Vertex{lower-\i}{\lowerMark}%
    \Vertex{upper-\i}{\upperMark}%
  }%
}

\newcommand{\LowerLabel}[2]{\node[vertex label,below=4pt] at (#1,0) {$#2$};}
\newcommand{\UpperLabel}[2]{\node[vertex label,above=4pt] at (#1,1) {$#2$};}
\newcommand{\CommonStart}[1]{%
  \node[case label] at (0,2) {$O_m,\qquad m\equiv #1\pmod 4,\qquad m>4$};
  \Comb{x/o,x/x,o/x,x/x,x/o,x/x,o/x,x/x}
  \LowerLabel{0}{v_1}
  \LowerLabel{1}{v_2}
  \LowerLabel{3}{v_4}
  \UpperLabel{1}{w_2}
  \UpperLabel{2}{w_3}
  \foreach \i in {7.5,8,8.5}{\fill (\i,0) circle[radius=0.7pt];}
}

\tikzset{
    domvertex/.style = {
        circle,
        draw = black,
        fill = black,
        minimum size = 7pt,
        inner sep = 0pt
    },
    vertex/.style = {
        circle,
        draw = black,
        fill = white,
        minimum size = 7pt,
        inner sep = 0pt
    },
    edge/.style = {
        line width = 0.5pt
    }
}

\newtheorem{thm}{Theorem}[section]

\newtheorem{cor}[thm]{Corollary}
\newtheorem{prop}[thm]{Proposition}
\newtheorem{lemma}[thm]{Lemma}

\newtheorem{obs}[thm]{Observation}

\graphicspath{ {./figures/} }

\allowdisplaybreaks

\title{Minimal super dominating sets of maximum cardinality}

\author[1]{
    Diana Avella-Alaminos
    \thanks{Email: \texttt{avella@ciencias.unam.mx}}
}

\author[2]{
    Magdalena Lema\'nska
    \thanks{Email: \texttt{magleman@pg.edu.pl}}
}

\author[1]{
   Rita Zuazua
    \thanks{Email: \texttt{ritazuazua@ciencias.unam.mx}}
}

\affil[1]{Universidad Nacional Aut\'{o}noma de M\'{e}xico,
Ciudad de M\'{e}xico, M\'{e}xico}

\affil[2]{Faculty of Applied Physics and Mathematics, Gda\'nsk Tech, Poland}

\date{}

\begin{document}

\maketitle

\begin{abstract}
Let $G = (V,E)$ be a graph.  A subset $S \subseteq V$ of vertices is a {\it super dominating set} if every vertex $v \in V - S$ has a neighbour $u \in S$, whose only neighbour in $V - S$ is $v$.  A set $S \subseteq V$ is called {\it open irredundant} if every vertex $u \in S$ has at least one neighbour $v \in V - S$ whose only neighbour in $S$ is $u$. In this paper we establish connections between super dominating sets and open irredundant sets, similar to the connections that exist between dominating sets and irredundant sets. In particular,  we define the new parameter $\Gamma _{sp}(G)$ as the maximum cardinality of a minimal super dominating set of $G$ and  calculate $\Gamma _{sp}(G)$ for some families of graphs $G$. Finally we present some properties of open irredundant sets and we will develop Gallai theorems involving domination, super domination and open irredundance parameters of graphs.
\end{abstract}

\section{Introduction}

Let $G=(V,E)$ a simple graph. A set vertex set $S \subseteq V$ is called {\it independent} if no two vertices in $S$ are adjacent.  The {\it vertex independence number} $\alpha(G)$ equals the maximum cardinality of an independent set of $G$. A set $S \subseteq V$ is called a {\it vertex cover} if for every edge $uv \in E$, either $u$, or $v$, or both $u$ and $v$, are in $S$. The {\it vertex covering number} $\beta(G)$ equals the minimum cardinality of a vertex cover of $G$.
In 1959 Gallai \cite{Gallai59} observed that a set $S \subset V$ is a maximal independent set if and only if the complement $V - S$ is a minimal vertex cover.  From this it is straightforward to prove the following theorem.

\begin{thm} (Gallai \cite{Gallai59})
For any graph $G$ of order $n = |V|$, $\alpha(G) + \beta(G) = n$.
\end{thm}

Since 1959 many theorems of the form $A(G) + B(G) = n$ have been proven, where $A(G)$ and $B(G)$ are two maximum and minimum types of parameters of graphs.  Results of this form are now called {\it Gallai theorems}.  In 1988 Cockayne et al. \cite{Cockayne88} surveyed many Gallai theorems.  Even more Gallai theorems are found in a 2002 paper by Balasubramanian et al. \cite {Balasubramanian02}.

The {\it open neighbourhood} $N(u)$ of a vertex $u$ in a graph $G$ is the set of vertices $N(u) = \{v : uv \in E \}$ which are adjacent to $u$. The {\it closed neighbourhood} $N[u]$ is the set $N[u] = N(u) \cup \{u\}$. The {\it open (closed) neighbourhood} $N(S)$  ($N[S]$) of a set $S$ of vertices is the union of the open (closed) neighbourhoods $N(u)$ ($N[u]$) of vertices $u \in S$. 
For a set $S$ of vertices, we say that a vertex $u \in S$ has an {\it external private neighbour with respect to $S$} if $N(u) - N[S - \{u\}] \ne \emptyset$. 

In this paper we are interested in the following three types of set of vertices in a graph.

A set $S$ of vertices is a {\it dominating set} if $N[S] = V,$ or equivalently, if every vertex $v \in V - S$ has at least one neighbour in $S$. The {\it domination number} $\gamma(G)$ of a graph $G$ equals the minimum cardinality of a dominating set of $G$. A dominating set of $G$ with minimum cardinality is called a {\it $\gamma$-set} of $G$.

A set $S \subseteq V$ is called a {\it super dominating set of $G$} if every vertex $v \in V - S$ has an external private neighbour $u \in S$, with respect to $V - S$. The {\it super domination number} $\gamma_{sp}(G)$ equals the minimum cardinality of a super dominating set in a graph $G$. A super dominating set of $G$ with minimum cardinality is called a {\it $\gamma_{sp}$-set} of $G$.

A set $S \subseteq V$ is called a {\it open irredundant set of $G$} if
every vertex $u \in S$ has an external private neighbour with respect to $S$. The {\it open irredundance number} $oir(G)$ equals the minimum cardinality of a maximal open irredundant set of $G$. An open irredundant set of $G$ with minimum cardinality is called an {\it $oir$-set} of $G$.  Open irredundance was considered in \cite{Cockayne03, Cockayne04, Cockayne08, Farley84,Farley83}.

The paper is organized in the following way: in the second section we recall the basic properties of a super domination in $G$ and we define the new parameter $\Gamma _{sp}(G)$ as the maximum cardinality of a minimal super dominating set of $G$. We calculate $\Gamma _{sp}(G)$ for some families of graphs $G.$ 
In section three we present some properties of open irredundant set and we will develop Gallai theorems involving domination, super domination and open irredundance parameters of graphs.

\section{Minimal super dominating sets of maximum cardinality}

A set $S \subseteq V$ is a {\it super dominating set of $G$} if for every vertex $v \in V - S$, there exists a vertex $u \in S$ such that $N(u) \cap (V - S) = \{v\}$. In the other words, $S$ is a super dominating set if every vertex $v \in V - S$ has an external private neighbour $u \in S$, with respect to $V - S$. The {\it super domination number} $\gamma_{sp}(G)$ equals the minimum cardinality of a super dominating set in a graph $G$.  
 Super dominating sets were defined in 2015 by Lema\'{n}ska, Swaminathan, Venkatakrishnan and Zuazua \cite{Lemanska15}. In the paper the authors calculate the super dominating number for cycles and paths and observe that for any graph $G$ with order $n$, $\gamma _{sp}(G)\geq \frac{n}{2}$. In Theorem 2.1, they give a characterization for connected graph $G$ with  $\gamma _{sp}(G)=\frac{n}{2}$. Moreover, in section 3, they show a characterization of trees which attain the lower and the upper bounds in terms of $n$. 

More properties of super dominating sets are studied in  \cite{Detlaff19} and \cite{Klein20}.  Recent papers on super domination include those by Krishnamurti and Venkatakrishnan \cite{Krishnamurti16} in 2016, Rodr\'{i}guez-Vel\'{a}zquez, Klein, and Yi \cite{Rodriguez20} in 2020, Senthilkumar, Venkatakrishnan and Naresh Kumar \cite{Senthilkumar22} and Zhuang \cite{Zhuang22} in 2022, Bujt\'{a}s, Ghanbari and Klav\v{z}ar \cite{Bujtas23} and Cabrera-Mart\'{i}nez and Conchalo Peir\'{o}  \cite{Cabrera23} in 2023, and Ghanbari \cite{Ghanbari24}  and Ghanbari, J\"{a}ger and Lehtil\"{a} \cite{Ghanbari24b} in 2024.

We define the new parameter {\it upper super domination number} $\Gamma_{sp}(G)$ as the maximum cardinality of a minimal super dominating set of $G$ and a minimal superdominating set of $G$ with maximum cardinality  we call a $\Gamma_{sp}$-set of $G.$ By definition, for any graph $G$, $\gamma_{sp}(G) \le  \Gamma_{sp}(G)$.

In the next proposition we show that the difference between $\Gamma_{sp}(G)$ and $\gamma_{sp}(G)$ can be arbitrarily large.

\begin{prop}\label{propGammasp}
For any $k \in \mathbb{N}^{+}$ there exists a graph $G$ such that $\Gamma_{sp}(G) - \gamma_{sp}(G) = k$.
\end{prop}

\begin{proof}  Consider a graph $G$  from  Figure~1. On the left side there is $G$ with depicted in black vertices belonging to a $\gamma_{sp}(G)$-set and in the left side $G$ with depicted in black vertices belonging to a $\Gamma_{sp}(G)$-set.

\begin{figure}[h]\label{rys1}
\begin{center}
\begin{tikzpicture}[scale=0.6]

\draw(0,0) -- (1,0)--(2,-0.5)--(3,-0.5)--(4,0)--(5,0)--(5,-1)--(4,-1)--(3,-0.5)--(2,-0.5)--(1,-1)--(0,-1)--(0,0);
\draw(1,-1)--(1,-2);
\draw(4,-1)--(4,-2);

\filldraw[fill opacity=0.9,fill=black]  (0,0) circle (0.1cm);
\filldraw[fill opacity=0.9,fill=white]  (1,0) circle (0.1cm);
\filldraw[fill opacity=0.9,fill=white]  (2,-0.5) circle (0.1cm);
\filldraw[fill opacity=0.9,fill=black]  (3,-0.5) circle (0.1cm);
\filldraw[fill opacity=0.9,fill=black]  (4,0) circle (0.1cm);
\filldraw[fill opacity=0.9,fill=white]  (5,0) circle (0.1cm);

\filldraw[fill opacity=0.9,fill=black]  (0,-1) circle (0.1cm);
\filldraw[fill opacity=0.9,fill=white]  (1,-1) circle (0.1cm);
\filldraw[fill opacity=0.9,fill=black]  (4,-1) circle (0.1cm);
\filldraw[fill opacity=0.9,fill=white]  (5,-1) circle (0.1cm);
\draw(1,-0.1)--(1,-0.9);
\draw(4,0)--(4,-1);

\draw(1,-1.9)--(1,-2.9);
\draw(4,-2)--(4,-3);

\draw(1,-4.5)--(1,-5.5);
\draw(4,-4.5)--(4,-5.5);

\draw(0,-2) -- (1,-2)--(2,-2.5)--(3,-2.5)--(4,-2)--(5,-2)--(5,-3)--(4,-3)--(3,-2.5)--(2,-2.5)--(1,-3)--(0,-3)--(0,-2);

\filldraw[fill opacity=0.9,fill=black]  (0,-2) circle (0.1cm);
\filldraw[fill opacity=0.9,fill=white]  (1,-2) circle (0.1cm);
\filldraw[fill opacity=0.9,fill=white]  (2,-2.5) circle (0.1cm);
\filldraw[fill opacity=0.9,fill=black]  (3,-2.5) circle (0.1cm);
\filldraw[fill opacity=0.9,fill=black]  (4,-2) circle (0.1cm);
\filldraw[fill opacity=0.9,fill=white]  (5,-2) circle (0.1cm);

\filldraw[fill opacity=0.9,fill=black]  (0,-3) circle (0.1cm);
\filldraw[fill opacity=0.9,fill=white]  (1,-3) circle (0.1cm);
\filldraw[fill opacity=0.9,fill=black]  (4,-3) circle (0.1cm);
\filldraw[fill opacity=0.9,fill=white]  (5,-3) circle (0.1cm);

\node at (1,-3.5) {$\cdot$};
\node at (1,-3.7) {$\cdot$};
\node at (1,-3.9) {$\cdot$};

\node at (4,-3.5) {$\cdot$};
\node at (4,-3.7) {$\cdot$};
\node at (4,-3.9) {$\cdot$};

\node at (-0.50,-0.5) {$1$};
\node at (-0.50,-2.5) {$2$};
\node at (-0.50,-5) {$k$};

\draw(0,-4.5) -- (1,-4.5)--(2,-5)--(3,-5)--(4,-4.5)--(5,-4.5)--(5,-5.5)--(4,-5.5)--(3,-5)--(2,-5)--(1,-5.5)--(0,-5.5)--(0,-5.5);
\draw(0,-4.5)--(0,-5.5);

\filldraw[fill opacity=0.9,fill=black]  (0,-4.5) circle (0.1cm);
\filldraw[fill opacity=0.9,fill=white]  (1,-4.5) circle (0.1cm);
\filldraw[fill opacity=0.9,fill=white]  (2,-5) circle (0.1cm);
\filldraw[fill opacity=0.9,fill=black]  (3,-5) circle (0.1cm);
\filldraw[fill opacity=0.9,fill=black]  (4,-4.5) circle (0.1cm);
\filldraw[fill opacity=0.9,fill=white]  (5,-4.5) circle (0.1cm);

\filldraw[fill opacity=0.9,fill=black]  (0,-5.5) circle (0.1cm);
\filldraw[fill opacity=0.9,fill=white]  (1,-5.5) circle (0.1cm);
\filldraw[fill opacity=0.9,fill=black]  (4,-5.5) circle (0.1cm);
\filldraw[fill opacity=0.9,fill=white]  (5,-5.5) circle (0.1cm);

\draw(6,0) -- (7,0)--(8,-0.5)--(9,-0.5)--(10,0)--(11,0)--(11,-1)--(10,-1)--(9,-0.5)--(8,-0.5)--(7,-1)--(6,-1)--(6,0);
\draw(7,-1)--(7,-2);
\draw(10,-1)--(10,-2);

\filldraw[fill opacity=0.9,fill=white]  (6,0) circle (0.1cm);
\filldraw[fill opacity=0.9,fill=black]  (7,0) circle (0.1cm);
\filldraw[fill opacity=0.9,fill=black]  (8,-0.5) circle (0.1cm);
\filldraw[fill opacity=0.9,fill=black]  (9,-0.5) circle (0.1cm);
\filldraw[fill opacity=0.9,fill=black]  (10,0) circle (0.1cm);
\filldraw[fill opacity=0.9,fill=white]  (11,0) circle (0.1cm);

\filldraw[fill opacity=0.9,fill=white]  (6,-1) circle (0.1cm);
\filldraw[fill opacity=0.9,fill=black]  (7,-1) circle (0.1cm);
\filldraw[fill opacity=0.9,fill=black]  (10,-1) circle (0.1cm);
\filldraw[fill opacity=0.9,fill=white]  (11,-1) circle (0.1cm);
\draw(7,-0.1)--(7,-0.9);
\draw(10,0)--(10,-0.9);

\draw(6,-2) -- (7,-2)--(8,-2.5)--(9,-2.5)--(10,-2)--(11,-2)--(11,-3)--(10,-3)--(9,-2.5)--(8,-2.5)--(7,-3)--(6,-3)--(6,-2);

\filldraw[fill opacity=0.9,fill=white]  (6,-2) circle (0.1cm);
\filldraw[fill opacity=0.9,fill=black]  (7,-2) circle (0.1cm);
\filldraw[fill opacity=0.9,fill=black]  (8,-2.5) circle (0.1cm);
\filldraw[fill opacity=0.9,fill=black]  (9,-2.5) circle (0.1cm);
\filldraw[fill opacity=0.9,fill=black]  (10,-2) circle (0.1cm);
\filldraw[fill opacity=0.9,fill=white]  (11,-2) circle (0.1cm);

\filldraw[fill opacity=0.9,fill=white]  (6,-3) circle (0.1cm);
\filldraw[fill opacity=0.9,fill=black]  (7,-3) circle (0.1cm);
\filldraw[fill opacity=0.9,fill=black]  (10,-3) circle (0.1cm);
\filldraw[fill opacity=0.9,fill=white]  (11,-3) circle (0.1cm);

\node at (7,-3.5) {$\cdot$};
\node at (7,-3.7) {$\cdot$};
\node at (7,-3.9) {$\cdot$};

\node at (10,-3.5) {$\cdot$};
\node at (10,-3.7) {$\cdot$};
\node at (10,-3.9) {$\cdot$};
\draw(7,-2.1)--(7,-2.9);
\draw(10,-2.1)--(10,-2.9);

\draw(6,-4.5) -- (7,-4.5)--(8,-5)--(9,-5)--(10,-4.5)--(11,-4.5)--(11,-5.5)--(10,-5.5)--(9,-5)--(8,-5)--(7,-5.5)--(6,-5.5)--(6,-4.5);

\filldraw[fill opacity=0.9,fill=white]  (6,-4.5) circle (0.1cm);
\filldraw[fill opacity=0.9,fill=black]  (7,-4.5) circle (0.1cm);
\filldraw[fill opacity=0.9,fill=black]  (8,-5) circle (0.1cm);
\filldraw[fill opacity=0.9,fill=black]  (9,-5) circle (0.1cm);
\filldraw[fill opacity=0.9,fill=black]  (10,-4.5) circle (0.1cm);
\filldraw[fill opacity=0.9,fill=white]  (11,-4.5) circle (0.1cm);

\filldraw[fill opacity=0.9,fill=white]  (6,-5.5) circle (0.1cm);
\filldraw[fill opacity=0.9,fill=black]  (7,-5.5) circle (0.1cm);
\filldraw[fill opacity=0.9,fill=black]  (10,-5.5) circle (0.1cm);
\filldraw[fill opacity=0.9,fill=white]  (11,-5.5) circle (0.1cm);
\draw(7,-4.5)--(7,-5.5);
\draw(10,-4.5)--(10,-5.5);

\node at (2,-6.5) {$\gamma_{sp}(G)=5k$};
\node at (8,-6.5) {$\Gamma_{sp}(G)=6k$};
\end{tikzpicture}
\end{center}
\caption{Graph $G$ for which $\Gamma_{sp}(G)-\gamma_{sp}(G)=k.$ }
\end{figure}
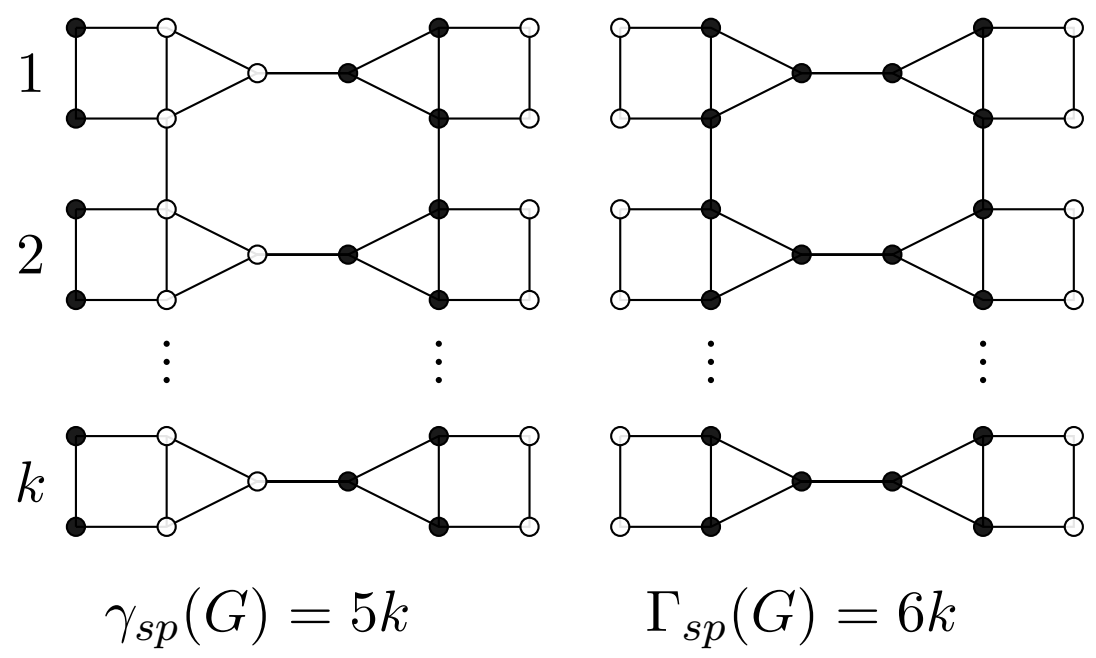
\end{proof}

We now determine the upper super domination number of the path on $n$ vertices. We begin with determining this parameter for the first $21$ paths.
The following table show the exact value of $\Gamma  _{sp}(P_n)$ for $2\leq n\leq 21.$

\begin{table}[ht]
\begin{center}
\resizebox{0.5\linewidth}{!}{\begin{tabular}{| c | c | c | c |c|c|c|c|c|c|c| }
\hline
$n$ & 2&3&4&5&6&7&8&9&10&11  \\ \hline
$\Gamma  _{sp}(P_n)$& 1&2&2&3&4&5&5&6&7&7 \\ \hline \hline
$n$ & 12&13&14&15&16&17&18&19&20&21  \\ \hline
$\Gamma  _{sp}(P_n)$ & 8&9&10&10&11&12&12&13&14&15  \\ \hline
\end{tabular}}
\caption{Exact values of $\Gamma  _{sp}(P_n)$ for $2\leq n\leq 21.$}

\end{center}
\end{table}

The vertices on the picture below on the left side of the paths  form a minimal super dominating sets of maximum cardinality.

\vspace{0.4cm}

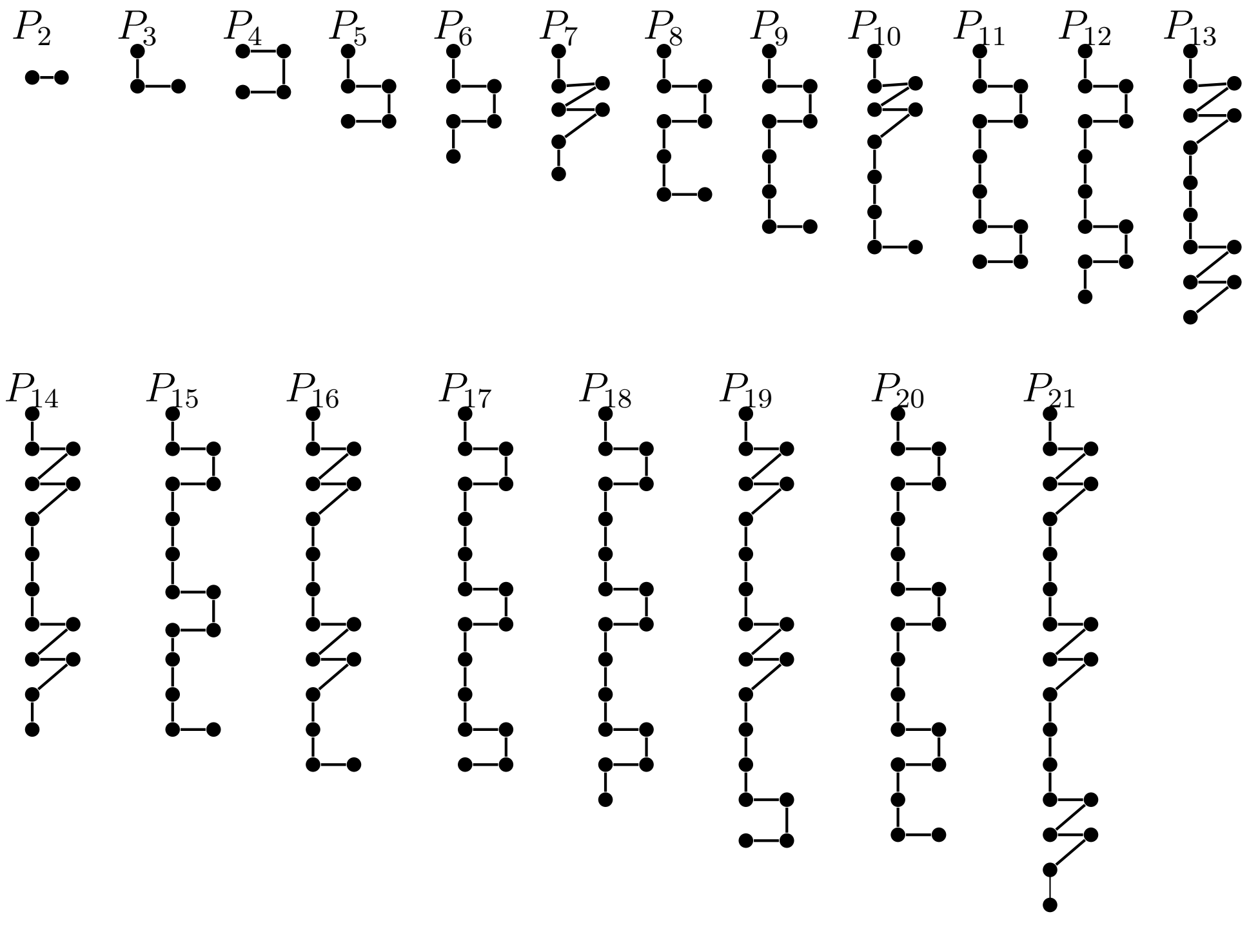
\tikzset{
  vertex/.style={circle, fill=black, inner sep=1.5pt},
  edge/.style={line width=0.8pt},
  num/.style={font=\large}
}

\newcommand{\Graph}[5]{
\begin{scope}[shift={(#2,#3)}]
  \node[num] at (0,0.85) {#1};

  \foreach \x/\y/\i in {#5}{
    \node[vertex] (G#1-\i) at (\x,\y) {};
  }

  \pgfmathtruncatemacro{\last}{#4-1}
  \foreach \i in {1,...,\last}{
    \pgfmathtruncatemacro{\j}{\i+1}
    \draw[edge] (G#1-\i)--(G#1-\j);
  }
\end{scope}
}

\begin{tikzpicture}[scale=0.6]

\Graph{$P_2$}{0}{0}{2}{
0/0/1, 0.5/0/2}

\Graph{$P_3$}{1.8}{0}{3}{
0/0.45/1, 0/-0.15/2, 0.7/-0.15/3}

\Graph{$P_4$}{3.6}{0}{4}{
 0/0.45/1, 0.7/0.45/2, 0.7/-0.25/3, 0/-0.25/4}

\Graph{$P_5$}{5.4}{0}{5}{
0/0.45/1, 0/-0.15/2, 0.7/-0.15/3, 0.7/-0.75/4, 0/-0.75/5}

\Graph{$P_6$}{7.2}{0}{6}{
0/0.45/1, 0/-0.15/2, 0.7/-0.15/3, 0.7/-0.75/4, 0/-0.75/5, 0/-1.35/6}

\Graph{$P_7$}{9.0}{0}{7}{
0/0.45/1, 0/-0.15/2, 0.75/-0.1/3, 0/-0.55/4, 0.75/-0.55/5, 0/-1.10/6, 0/-1.65/7}

\Graph{$P_8$}{10.8}{0}{8}{
0/0.45/1, 0/-0.15/2, 0.7/-0.15/3, 0.7/-0.75/4, 0/-0.75/5, 0/-1.35/6, 0/-2/7, 0.7/-2/8}

\Graph{$P_9$}{12.6}{0}{9}{
0/0.45/1, 0/-0.15/2, 0.7/-0.15/3, 0.7/-0.75/4, 0/-0.75/5, 0/-1.35/6, 0/-1.95/7, 0/-2.55/8, 0.7/-2.55/9}

\Graph{$P_{10}$}{14.4}{0}{10}{
0/0.45/1, 0/-0.15/2, 0.7/-0.1/3, 0/-0.55/4, 0.7/-0.55/5, 0/-1.10/6, 0/-1.70/7, 0/-2.30/8, 0/-2.90/9, 0.7/-2.90/10}

\Graph{$P_{11}$}{16.2}{0}{11}{
0/0.45/1, 0/-0.15/2, 0.7/-0.15/3, 0.7/-0.75/4, 0/-0.75/5, 0/-1.35/6, 0/-1.95/7, 0/-2.55/8, 0.7/-2.55/9, 0.7/-3.15/10, 0/-3.15/11}

\Graph{$P_{12}$}{18.0}{0}{12}{
0/0.45/1, 0/-0.15/2, 0.7/-0.15/3, 0.7/-0.75/4, 0/-0.75/5, 0/-1.35/6, 0/-1.95/7, 0/-2.55/8, 0.7/-2.55/9, 0.7/-3.15/10, 0/-3.15/11, 0/-3.75/12}

\Graph{$P_{13}$}{19.8}{0}{13}{
0/0.45/1, 0/-0.15/2, 0.75/-0.1/3, 0/-0.65/4, 0.75/-0.65/5, 0/-1.2/6, 0/-1.8/7, 0/-2.35/8, 0/-2.9/9, 0.75/-2.9/10, 0/-3.5/11, 0.75/-3.5/12, 0/-4.1/13}

\Graph{$P_{14}$}{0}{-6.2}{14}{
0/0.45/1, 0/-0.15/2, 0.7/-0.15/3, 0/-0.75/4, 0.7/-0.75/5,
0/-1.35/6, 0/-1.95/7, 0/-2.55/8, 0/-3.15/9, 0.7/-3.15/10,
0/-3.75/11, 0.7/-3.75/12, 0/-4.35/13, 0/-4.95/14}

\Graph{$P_{15}$}{2.4}{-6.2}{15}{
0/0.45/1, 0/-0.15/2, 0.7/-0.15/3, 0.7/-0.75/4, 0/-0.75/5,
0/-1.35/6, 0/-1.95/7, 0/-2.6/8, 0.7/-2.6/9, 0.7/-3.25/10,
0/-3.25/11, 0/-3.75/12, 0/-4.35/13, 0/-4.95/14, 0.7/-4.95/15}

\Graph{$P_{16}$}{4.8}{-6.2}{16}{
0/0.45/1, 0/-0.15/2, 0.7/-0.15/3, 0/-0.75/4, 0.7/-0.75/5,
0/-1.35/6, 0/-1.95/7, 0/-2.55/8, 0/-3.15/9, 0.7/-3.15/10,
0/-3.75/11, 0.7/-3.75/12, 0/-4.35/13, 0/-4.95/14, 0/-5.55/15, 0.7/-5.55/16}

\Graph{$P_{17}$}{7.4}{-6.2}{17}{
0/0.45/1, 0/-0.15/2, 0.7/-0.15/3, 0.7/-0.75/4, 0/-0.75/5,
0/-1.35/6, 0/-1.95/7, 0/-2.55/8, 0.7/-2.55/9, 0.7/-3.15/10,
0/-3.15/11, 0/-3.75/12, 0/-4.35/13, 0/-4.95/14, 0.7/-4.95/15,
0.7/-5.55/16, 0/-5.55/17}

\Graph{$P_{18}$}{9.8}{-6.2}{18}{
0/0.45/1, 0/-0.15/2, 0.7/-0.15/3, 0.7/-0.75/4, 0/-0.75/5,
0/-1.35/6, 0/-1.95/7, 0/-2.55/8, 0.7/-2.55/9, 0.7/-3.15/10,
0/-3.15/11, 0/-3.75/12, 0/-4.35/13, 0/-4.95/14, 0.7/-4.95/15,
0.7/-5.55/16, 0/-5.55/17, 0/-6.15/18}

\Graph{$P_{19}$}{12.2}{-6.2}{19}{
0/0.45/1, 0/-0.15/2, 0.7/-0.15/3, 0/-0.75/4, 0.7/-0.75/5,
0/-1.35/6, 0/-1.95/7, 0/-2.55/8, 0/-3.15/9, 0.7/-3.15/10,
0/-3.75/11, 0.7/-3.75/12, 0/-4.35/13, 0/-4.95/14, 0/-5.55/15,
0/-6.15/16, 0.7/-6.15/17, 0.7/-6.85/18, 0/-6.85/19}

\Graph{$P_{20}$}{14.8}{-6.2}{20}{
0/0.45/1, 0/-0.15/2, 0.7/-0.15/3, 0.7/-0.75/4, 0/-0.75/5,
0/-1.35/6, 0/-1.95/7, 0/-2.55/8, 0.7/-2.55/9, 0.7/-3.15/10,
0/-3.15/11, 0/-3.75/12, 0/-4.35/13, 0/-4.95/14, 0.7/-4.95/15,
0.7/-5.55/16, 0/-5.55/17, 0/-6.15/18, 0/-6.75/19, 0.7/-6.75/20}

\Graph{$P_{21}$}{17.4}{-6.2}{20}{
0/0.45/1, 0/-0.15/2, 0.7/-0.15/3, 0/-0.75/4, 0.7/-0.75/5,
0/-1.35/6, 0/-1.95/7, 0/-2.55/8, 0/-3.15/9, 0.7/-3.15/10,
0/-3.75/11, 0.7/-3.75/12, 0/-4.35/13, 0/-4.95/14, 0/-5.55/15,
0/-6.15/16, 0.7/-6.15/17, 0/-6.75/18, 0.7/-6.75/19, 0/-7.35/20, 0/-7.95/21}
\draw (17.4,-13.55) -- (17.4,-14.15);
\end{tikzpicture}

\vspace{0.3cm}

For $n\in\mathbb{N}^{+}$, let $P_{n}$ be a path with set of vertices $\{v_1,\dots,v_n\}$ and set of edges $\{\{v_i,v_{i+1}\}\,|\,i=1,\dots,n-1\}.$ For $i,j\in\{1,\dots,n-1\}$ with $i<j$, we will denote the subset $\{v_i,v_{i+1},\dots , v_{j}\}$ by $[v_i,v_j].$

\begin{obs}\label{obs:Pn1}
Let $n\in\mathbb{N}^{+}$, $n>2$. For every minimal super dominating set $S$ of $P_n$, $|S\cap [v_{n-2},v_{n}]| \leq 2.$
\end{obs}
\begin{proof}
If $|S\cap  [v_{n-2},v_{n}]| =3,$ then $S\setminus\{v_{n}\}$ would be a super dominating set of $P_n$ contradicting that $S$ is minimal.
\end{proof}

\begin{obs}\label{obs:Pn2}
Let $n\in\mathbb{N}^{+}$, $n>4$. For every minimal superdominating set $S$ of $P_n$ and every $i\in\{1,\dots,n-4\}$, $|S\cap [v_{i},v_{i+4}]| \leq 4.$
\end{obs}
\begin{proof}
If $|S\cap  [v_{i},v_{i+4}]| =5,$ then $S\setminus\{v_{i+2}\}$ would be a super dominating set of $P_n$ contradicting that $S$ is minimal.
\end{proof}

\begin{obs}\label{obs:Pn3}
Let $n\in\mathbb{N}^{+}$, $n>5$. For every minimal superdominating set $S$ of $P_n$, $|S\cap  [v_{n-5},v_{n}]| \leq 4.$
\end{obs}
\begin{proof}
By Observation \ref{obs:Pn2}, $|S\cap [v_{n-5},v_{n}]| <6$. Suppose $|S\cap [v_{n-5},v_{n}]| =5$, by Observation \ref{obs:Pn1}, $v_{n-5},v_{n-4},v_{n-3}\in S$, and by Observation \ref{obs:Pn2}, $v_{n}\in S$. If $v_{n-2}\notin S$, then $S\setminus\{v_{n-3}\}$ would be a superdominating set of $P_n$ contradicting that $S$ is minimal, and if $v_{n-1}\notin S$, then $S\setminus\{v_{n-2}\}$ would be a super dominating set of $P_n$ contradicting that $S$ is minimal.
\end{proof}

\begin{lemma}\label{lem:Pn}
Let $n\in\mathbb{N}^{+}$, $n>6$. For every minimal superdominating set $S$ of $P_n$ and every $i\in\{1,\dots,m-6\}$, $|S\cap [v_{i},v_{i+6}]| \leq 5.$
\end{lemma}
\begin{proof}
Suppose $|S\cap [v_{i},v_{i+6}]| =7,$ then $S\setminus\{v_{i+3}\}$ would be a super dominating set of $P_n$ contradicting that $S$ is minimal.\\ Now, suppose $|S\cap [v_{i},v_{i+6}]| =6$, by Lemma \ref{obs:Pn2},  none five consecutive vertices belong to $S$, so $v_i,v_{i+1},v_{i+2},v_{i+5}\,v_{i+6}\in S$, if $v_{i+2}\notin S$ or $v_{i+4}\notin S$, then $S\setminus\{v_{i+3}\}$ would be a super dominating set of $P_n$, and if $v_{i+3}\notin S$, then $S\setminus\{v_{i+4}\}$ would be a super dominating set of $P_n$, contradicting in both cases that $S$ is minimal.
\end{proof}

\begin{prop}\label{prop:Pm1}
Let $n\in\mathbb{N}^{+}$. Let $q,r\in \mathbb{N}$ such that $n=7q+r$, with $0\leq r<7$. Then $\Gamma_{sp}(P_n)\leq 5q+\lceil\frac{2r}{3}\rceil .$
\end{prop}
\begin{proof}
Let $S$ be a minimal super dominating set of $P_n$. We consider some cases, depending from the value of $r$.

\vspace{0.2cm}

{\bf Case 1:  $r=0$.}  By Lemma \ref{lem:Pn}, every $i\in\{1,8,\dots,7(q-1)+1\}$, $|S\cap [v_{i},v_{i+6}]| \leq 5,$ then 
$$|S|=\sum_1^{7(q-1)+1}|S\cap[v_{i},v_{i+6}]| \leq 5q=5q+\left\lceil\frac{2(0)}{3}\right\rceil.$$

For the remaining cases, we see that if $q=0$ the result is clear, because  $\Gamma_{sp}(P_1)=\Gamma_{sp}(P_2)=1,\Gamma_{sp}(P_3)=\Gamma_{sp}(P_4)=2, \Gamma_{sp}(P_5)=3,\Gamma_{sp}(P_6)=4$ and all those values are less or equal to $\lceil\frac{2r}{3}\rceil$. Suppose now $q>0$.

\vspace{0.2cm}

{\bf Case 2:  $r=1$.}  By Lemma \ref{lem:Pn}, for every $i\in\{1,8,\dots,7(q-1)+1\}$, $|S\cap [v_{i},v_{i+6}]| \leq 5,$ then 
$$|S|=\left( \sum_1^{7(q-1)+1}|S\cap[v_{i},v_{i+6}]|\right) +| S\cap\{v_n\}| \leq 5q+1=5q+\left\lceil\frac{2(1)}{3}\right\rceil.$$

{\bf Case 3:  $r=2$.}  By Lemma \ref{lem:Pn}, for every $i\in\{1,8,\dots,7(q-1)+1\}$, $|S\cap [v_{i},v_{i+6}]| \leq 5,$ then 
$$|S|=\left( \sum_1^{7(q-1)+1}|S\cap[v_{i},v_{i+6}]| \right)+| S\cap\{v_{n-1},v_n\}| \leq 5q+2=5q+\left\lceil\frac{2(2)}{3}\right\rceil.$$

{\bf Case 4:  $r=3$.}  By Lemma \ref{lem:Pn}, for every $i\in\{1,8,\dots,7(q-1)+1\}$, $|S\cap [v_{i},v_{i+6}]| \leq 5,$ and by Observation \ref{obs:Pn1}, $| S\cap[v_{n-2},v_{n}]| \leq 2$, then 
$$|S|=\left( \sum_1^{7(q-1)+1}|S\cap[v_{i},v_{i+6}]| \right)+| S\cap[v_{n-2},v_{n}]| \leq 5q+2=5q+\left\lceil\frac{2(3)}{3}\right\rceil.$$

{\bf Case 5: $r=4$. } By Lemma \ref{lem:Pn}, for every $i\in\{1,8,\dots,7(q-1)+1\}$, $|S\cap [v_{i},v_{i+6}]| \leq 5,$ and by Observation \ref{obs:Pn1}, $| S\cap[v_{n-2},v_{n}]| \leq 2$, then 
$$\begin{array}{ll}
|S|
&=\left( \sum_1^{7(q-1)+1}|S\cap[v_{i},v_{i+6}]|\right) +| S\cap\{v_{n-3}\}|+| S\cap[v_{n-2},v_{n}]| \\&\leq 5q+1+2=5q+3=5q+\left\lceil\frac{2(4)}{3}\right\rceil.
\end{array}
$$

{\bf Case 6: $r=5$.}  By Lemma \ref{lem:Pn}, for every $i\in\{1,8,\dots,7(q-1)+1\}$, $|S\cap [v_{i},v_{i+6}]| \leq 5,$ and by Observation \ref{obs:Pn2}, $|S\cap [v_{n-4},v_{n}]| \leq 4$, then 
$$
|S|=\left( \sum_1^{7(q-1)+1}|S\cap[v_{i},v_{i+6}]|\right) +|S\cap [v_{n-4},v_{n}]| \leq 5q+4=5q+\left\lceil\frac{2(5)}{3}\right\rceil.
$$

{\bf Case 7: $r=6$.}  By Lemma \ref{lem:Pn}, for every $i\in\{1,8,\dots,7(q-1)+1\}$, $|S\cap [v_{i},v_{i+6}]| \leq 5,$ and by Observation \ref{obs:Pn3}, $|S\cap [v_{n-5},v_{n}]|  \leq 4$, then 
$$
|S|=\left( \sum_1^{7(q-1)+1}|S\cap[v_{i},v_{i+6}]| \right)+|S\cap [v_{n-5},v_{n}]| \leq 5q+4=5q+\left\lceil\frac{2(6)}{3}\right\rceil.
$$
In all cases $\Gamma_{sp}(P_n)\leq |S |\leq 5q+\lceil\frac{2r}{3}\rceil .$
\end{proof}

\begin{prop}\label{prop:Pm2}
Let $n\in\mathbb{N}^{+}$, $n\geq 5$. Let $q,r\in \mathbb{N}$ such that $n=21q+r$, with $0\leq r<21$. Then $\Gamma_{sp}(P_n)\geq 15q$ if $r=0,1,$ and $\Gamma_{sp}(P_n)\geq 15q+\lceil\frac{2r-1}{3}\rceil $ if $2\leq r<21.$
\end{prop}
\begin{proof}

For each $j\in \{0,1,\dots ,q\}$ let $B_0=\emptyset$ and \\
$B_j=\{v_{21(j-1)+k}\mid k\in\{1,2,4,6,7,8,9,11,13,14,15,16,18,20,21\}\}$ if $j>0$.

\vspace{0.2cm}

 If $r=0$ or $r=1,$ then the set $S=\cup_{j=0}^qB_j$ is a minimal super dominating set of $P_n$ with $|S|=15q$. Then, $\Gamma_{sp}(P_n)\geq |S|= 15q.$

\vspace{0.2cm}

 For $r\geq 2$ define: 
\begin{itemize}
\item if $r=16,$
$T_n=\{v_{21q+k}\mid k\in\{1,\dots r-1\},k\not\equiv 3(mod \, 7),k\not\equiv 5(mod \, 7)\},$
\item if $r\equiv 1(mod \, 3)$ and $r\neq 16,$
$T_n=\{v_{21q+k}\mid k\in\{1,\dots r\},k\not\equiv 3(mod \, 7),k\not\equiv 5(mod \, 7)\},$
\item if $r\equiv 2(mod \, 6)$, 
$T_n=\{v_{21q+k}\mid k\in\{1,\dots r-2\},k\not\equiv 3(mod \, 6),k\not\equiv 4(mod \, 6)\}\cup\{v_{n-1}\},$
\item if $r\not\equiv 1(mod \, 3)$ and $r\not\equiv 2(mod \, 6)$, 
$T_n=\{v_{21q+k}\mid k\in\{1,\dots r\},k\not\equiv 3(mod \, 6),k\not\equiv 4(mod \, 6)\}.$
\end{itemize}

The set $S=\left(\cup_{j=0}^qB_j\right)\cup T_n$ is a minimal super dominating set of $P_n$ with $|S|=15q+\lceil\frac{2r-1}{3}\rceil$. Then, $\Gamma_{sp}(P_n)\geq |S|=15q+\lceil\frac{2r-1}{3}\rceil .$
\end{proof}

\begin{thm}\label{prop:Pm3}
Let $n\in\mathbb{N}^{+}$. Let $q,r\in \mathbb{N}$ such that $n=21q+r$, with $0\leq r<20$ and $k,t\in \mathbb{N}$ such that $r=7k+t$, with $0\leq t<7$. Then $\Gamma_{sp}(P_n)= 15q+5k+\lceil\frac{2t}{3}\rceil$ or $\Gamma_{sp}(P_n)= 15q+5k+\lceil\frac{2t}{3}\rceil-1$.
\end{thm}
\begin{proof} We consider some cases, depending on the values of $r.$  The cases $n\in \{2,3,4\}$ are immediate, so we may assume that $n\geq 5$.

\vspace{0.2cm}

{\bf Case 1: $r=0.$}  In this case $n=21q+r$, with $r=0$, and $r=7k+t$, with $k=t=0$. As $n=21q+r=7(3q)+r$ with $r=0$, by Proposition \ref{prop:Pm1}, $\Gamma_{sp}(P_n)\leq 5(3q)+\left\lceil\frac{2r}{3}\right\rceil= 15q$, and by Proposition \ref{prop:Pm2}, $\Gamma_{sp}(P_n)\geq 15q$.  Then, $\Gamma_{sp}(P_n)= 15q=15q+5k+\lceil\frac{2t}{3}\rceil$. 

\vspace{0.2cm}

{\bf Case 2: $r=1.$}  In this case $n=21q+r$, with $r=1$, and $r=7k+t$, with $k=0$ and $t=1$. As $n=21q+r=7(3q)+r$ with $r=1$, by Proposition \ref{prop:Pm1}, $\Gamma_{sp}(P_n)\leq 5(3q)+\left\lceil\frac{2r}{3}\right\rceil= 15q+\left\lceil\frac{2}{3}\right\rceil$, and by Proposition \ref{prop:Pm2}, $\Gamma_{sp}(P_n)\geq 15q$, and $15q=15q+\left\lceil\frac{2}{3}\right\rceil-1$.  Then, $\Gamma_{sp}(P_n)= 15q+\lceil\frac{2}{3}\rceil=15q+\lceil\frac{2t}{3}\rceil$ or $\Gamma_{sp}(P_n)= 15q+\lceil\frac{2}{3}\rceil-1=15q+\lceil\frac{2t}{3}\rceil-1$. 

\vspace{0.2cm}

{\bf Case 3:  $2\leq r$.} 
As $n=21q+r$, with $2\leq r<20$ and $k,t\in \mathbb{N}$ such that $r=7k+t$, with $0\leq t<7$, then $n=7(3q)+7k+t=7(3q+k)+t$, so by Proposition \ref{prop:Pm1}, $$\Gamma_{sp}(P_n)\leq 5(3q+k)+\left\lceil\frac{2t}{3}\right\rceil =15q+5k+\left\lceil\frac{2t}{3}\right\rceil.$$
Now,  by Proposition \ref{prop:Pm2}, $$\begin{array}{ll}\Gamma_{sp}(P_n)&\geq 15q+\left\lceil\frac{2r-1}{3}\right\rceil=15q+\left\lceil\frac{2(7k+t)-1}{3}\right\rceil= 15q+\left\lceil\frac{15k+2t-k-1}{3}\right\rceil\\&= 15q+5k+\left\lceil\frac{2t-k-1}{3}\right\rceil,\end{array}$$
but, $r=7k+t$ and $2\leq r<20,$ so $k\in\{0,1,2\}$. Then, $\left\lceil\frac{2t-k-1}{3}\right\rceil=\left\lceil\frac{2t}{3}\right\rceil$, (when $k=0$ and $t\not\equiv 2(mod \, 3)$, or when $k=1$ and $t\equiv 0(mod \, 3)$), or $\left\lceil\frac{2t-k-1}{3}\right\rceil=\left\lceil\frac{2t}{3}\right\rceil -1$ in the remaining cases, so 
$$\Gamma_{sp}(P_n)\geq 15q+5k+\left\lceil\frac{2t}{3}\right\rceil$$ or $$\Gamma_{sp}(P_n)\geq15q+5k+\left\lceil\frac{2t}{3}\right\rceil-1.$$

We conclude that $\Gamma_{sp}(P_n)= 15q+5k+\lceil\frac{2t}{3}\rceil$ or $\Gamma_{sp}(P_n)= 15q+5k+\lceil\frac{2t}{3}\rceil-1$.

\end{proof}

By the analysis made in the proof of Theorem \ref{prop:Pm3}, we have the following

\begin{cor}
Let $n\in\mathbb{N}^{+}$. Let $q,r\in \mathbb{N}$ such that $n=21q+r$, with $0\leq r<20$ and $k,t\in \mathbb{N}$ such that $r=7k+t$, with $0\leq t<7$. In the following cases:
\begin{itemize}
\item[(a)] $r=0$, or
\item[(b)] $2\leq r$, $k=0$ and $t\not\equiv 2(mod \, 3)$, or
\item[(c)] $2\leq r$, $k=1$ and $t\equiv 0(mod \, 3)$,
\end{itemize}
we have that $\Gamma_{sp}(P_n)= 15q+5k+\lceil\frac{2t}{3}\rceil$. 
\end{cor}

Now we consider the case of cycles. 
For $n\in\mathbb{N}^{+}$, let $C_{n}$ be the cycle with set of vertices $\{v_1,\dots,v_n\}$ and set of edges $\{\{v_i,v_{i+1}\}\,|\,i=1,\dots,n-1\}\cup \{v_n,v_1\}.$ 

\begin{lemma}\label{obs:Cn}
Let $n\in\mathbb{N}^{+}$, $n>2$. For every minimal super dominating set $S$ of $C_n$, there exist two adjacent vertices  $u,v\in V(G)$ of $C_n$  such that  $u,v\in S.$
\end{lemma}
\begin{proof}
Let $S$ be a minimal super dominating set of $C_n$. Consider $w \in V - S$; as $S$ is a super dominating set, $w$ has an external private neighbour $u\in S$ with respect to $V - S$. Then, as $u$ is an external private neighbour with respect to $V - S$, $v\in N(u) - \{w\}$ is in $S$. 
\end{proof}
\begin{lemma}\label{lem211}
Let $S$ be a $\Gamma _{sp}(C_n)$-set and let  $u,v$ be adjacent vertices in $S.$  Then any longest path $P$ in $C_n$  with $V(P)\subseteq S$ and $u,v\in V(P)$ has   a length  $l$ satisfying $1\leq l\leq 3.$
\end{lemma}
\begin{proof} Suppose by contradiction that $P$ has length more than three; then there exists a path $P_5=x_1x_2x_3x_4x_5$ contained in $P.$ But in that case $S- \{x_3\}$ is a super dominating set of $C_n$ with cardinality $|S|-1,$ a contradiction with minimality of $S$. Thus $P=P_m$ with $m\in \{ 2,3,4\}.$
\end{proof}
\begin{prop}
For every integer $n \geq 3, $ we have $\Gamma_{sp}(C_n)\leq \Gamma_{sp}(P_n)$. 
\end{prop}  
\begin{proof}
Let $S$ be a $\Gamma _{sp}(C_n)$-set and consider two vertices $u,v$ adjacent in $S.$  Let $P_m$ be a longest path with $V(P_{m})\subseteq S$ and $u,v\in V(P_{m})$; from Lemma \ref{lem211}  $m\in \{ 2,3,4\}$. 
\begin{itemize}
    \item If $m=2$, then $S$ is a minimal super dominating set of $P_n$, where $P_n=C_n-\{uv\}$;
    \item If $m=3$ and $P_3=x_1x_2x_3$, then $S$ is a minimal super dominating set of $P_n$, where $P_n=C_n-\{x_1x_2\}$; 
    \item If $m=4$ and $P_4=x_1x_2x_3x_4$, then $S$ is a minimal super dominating set of $P_n$, where $P_n=C_n-\{x_2x_3\}$.
    \end{itemize}
    Then, $\Gamma_{sp}(C_n)=|S|\leq \Gamma_{sp}(P_n).$
    \end{proof}

Finally, we consider minimal super dominating sets of maximum cardinality in caterpillars.  For $m\in\mathbb{N}^{+}$ with $m>2$, let $O_{m}$ be a caterpillar graph with set of vertices $\{v_1,\dots,v_m\}\cup\{w_1,\dots ,w_m\}$ and set of edges $\{\{v_i,v_{i+i}\}\,|\,i=1,\dots,m-1\}\cup \{\{v_i,w_{i}\}\,|\,i=1,\dots,m\}.$ Notice that $|V(O_m)|=2m.$ 

\begin{lemma}\label{lem:Cm}
Let $m\in\mathbb{N}^{+}$ with $m>2$. For every minimal super dominating set of $O_m$, $|S\cap \{v_1,w_1\}|=|S\cap \{v_m,w_m\}|=1.$
\end{lemma}
\begin{proof}
As $w_1$ must be dominated by $S$, then  $S\cap \{v_1,w_1\}\neq \emptyset$. Suppose by contradiction that $S\cap \{v_1,w_1\}=\{v_1,w_1\},$ then, as $S-\{w_1\}$ is not a super dominating set, $b_1$ would not have a private neighbour with respect to $V-(S-\{w_1\})$, so $v_2\notin S$, but this implies, as $w_2$ must be dominated by $S$, $w_2\in S$. Then $S-\{v_1\}$ would be a super dominating set, which is a contradiction.

We conclude $|S\cap \{v_1,w_1\}|=1.$ Similarly, $|S\cap \{v_m,w_m\}|=1.$ 
\end{proof}

\begin{lemma}\label{lem:Cm2}
Let $m\in\mathbb{N}^{+}$ with $m>2$. For every minimal super dominating set of $O_m$ and for every $i\in\{1,2,\dots ,m-1\}$, $|S\cap \{v_i,w_i,v_{i+1},w_{i+1}\}|\leq 3.$
\end{lemma}
\begin{proof}
Suppose by contradiction that $|S\cap \{v_i,w_i,v_{i+1},w_{i+1}\}|=4.$ By Lemma \ref{lem:Cm}, $i\neq 1$ and  $i\neq n-1$.
As $S-\{w_i\}$ is not a super dominating set, $v_{i-1}\notin S$. Also, as $S-\{w_{i+1}\}$ is not a super dominating set, $v_{i+1}\notin S$. As consequence, $w_{i-1},w_{i+1}\in S$, otherwise they would not have a private neighbour with respect to $V-S$. Then, $S-\{w_i,w_{i+1}\}$ would be a super dominating set, which is a contradiction.

We conclude $ |S\cap \{v_i,w_i,v_{i+1},w_{i+1}\}|\leq 3.$
\end{proof}

\begin{thm}\label{Cm}
Let $m>2$. Then, 
$ \Gamma_{sp}(O_{m})
= \left\{
    \begin{array}{ll}
        3q+1, &if\, m=2q+1, q\geq 1,\\
        3q+2, &if\, m=2q+2, q\geq 1.
    \end{array}
    \right.$
\end{thm}

\begin{proof}

First we will prove that $ \Gamma_{sp}(O_{m})
\geq \left\{
    \begin{array}{ll}
         3q+1, &if\, m=2q+1, q\geq 1,\\
        3q+2, &if\, m=2q+2, q\geq 1.
    \end{array}
    \right.$

\vspace{.3cm}
We consider some cases, depending on the value of $m.$ Minimal super dominating sets of maximum cardinality of $O_m$ are marked on Figure~\ref{GammaspCm}.

\vspace{0.2cm}

{\bf Case 1:   $m=2q+1,$ with $q\geq 1.$}
\begin{enumerate}
\item For $m\equiv 1\,(mod\,4)$ the set $$S_{m}= \{v_j\,|\,1\leq j\leq m,j\not\equiv 3\,(mod\,4)\}\cup \{w_{j}\,|\, 1\leq j< m,j\not\equiv 1\,(mod\,4)\}$$ is a minimal super dominating set of $O_m$, so $\Gamma_{sp}(O_{m})\geq |S_{m}|=6(\frac{m-1}{4})+1=6(\frac{2q}{4})+1=3q+1$. 
\item For $m\equiv 3
\,(mod\,4)$, the set $$S_{m}= \{v_j\,|\,1\leq j\leq m,j\not\equiv 3\,(mod\,4)\}\cup \{w_{j}\,|\,1\leq j\leq m,j\not\equiv 1\,(mod\,4)\}$$ is a minimal super dominating set of $O_m$, so $\Gamma_{sp}(O_{m})\geq |S_{m}|=6(\frac{m-3}{4})+4=6(\frac{2q-2}{4})+4=3(q-1)+4=3q+1$. 
\end{enumerate}

\vspace{0.2cm}

{\bf Case 2:  $m=2q+2,$ with $ q\geq 1$.}
\begin{itemize}
\item[(a)] For $m\equiv 0\,(mod\,4)$, the set 
$$\begin{array}{ll}S_{m}&= \{v_j\,|\,1\leq j\leq m-4,j\not\equiv 3\,(mod\,4)\}\cup \{w_{j}\,|\,1\leq j\leq m-4,j\not\equiv 1\,(mod\,4)\}\\&\phantom{=}\cup \{v_{m-3},v_{m-2},w_{m-2},w_{m-1},w_m\}\end{array}$$ 
is a minimal super dominating set of $O_m$, so $\Gamma_{sp}(O_{m})\geq |S_{m}|=6(\frac{m-4}{4})+5=6(\frac{2q-2}{4})+5=3(q-1)+5=3q+2$. 
\item[(b)] For $m\equiv 2\,(mod\,4)$, the set $$S_{m}= \{v_j\,|\,1\leq j\leq m,j\not\equiv 3\,(mod\,4)\}\cup \{w_{j}\,|\,1\leq j\leq m-2,j\not\equiv 1\,(mod\,4)\}$$ is a minimal super dominating set of $O_m$, so $\Gamma_{sp}(O_{m})\geq |S_{m}|=6(\frac{m-2}{4})+2=6(\frac{2q}{4})+2=3q+2$. 
\end{itemize}

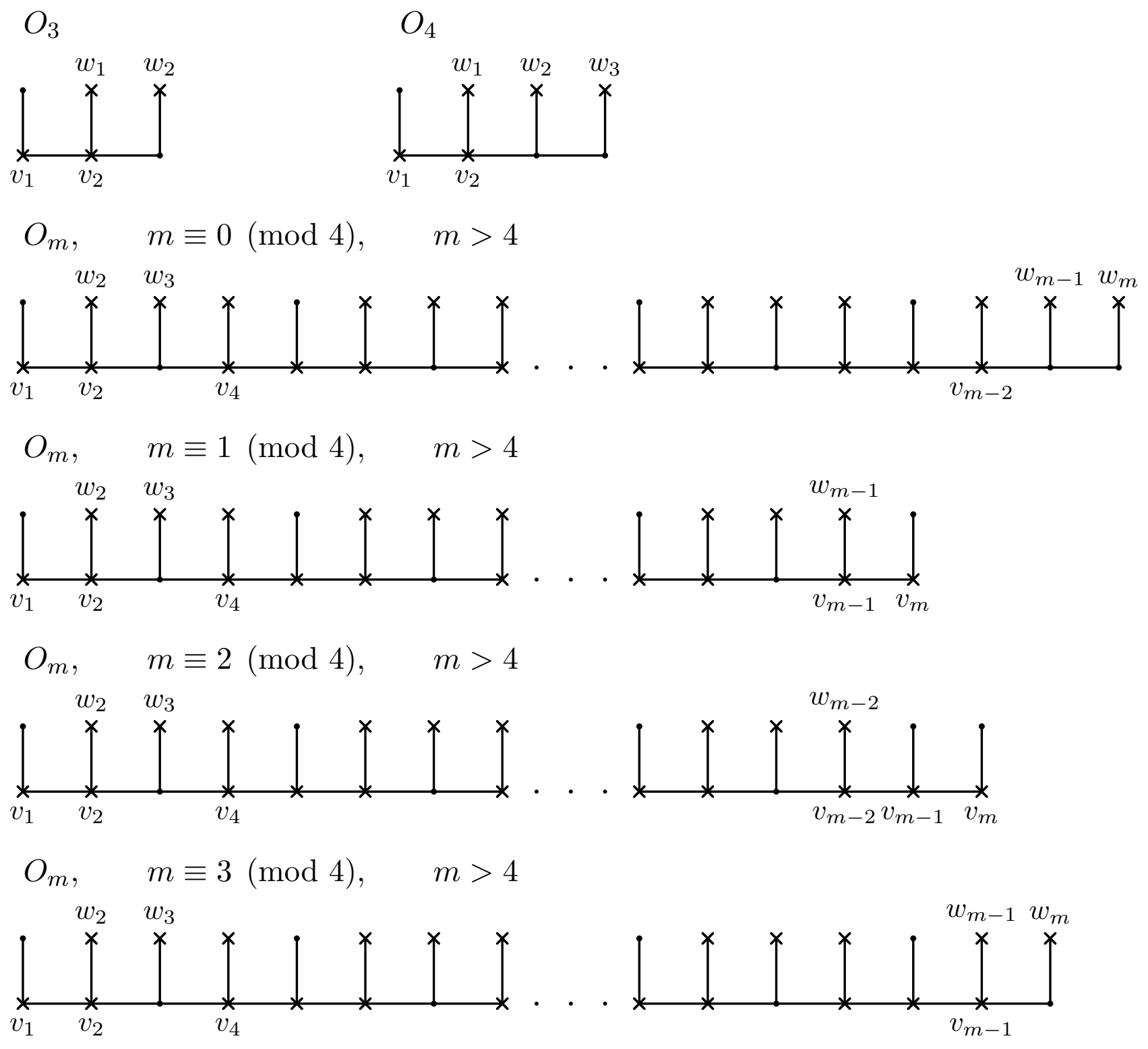
\begin{figure}[htbp]
\centering
\scalebox{0.8}{\begin{tikzpicture}[x=8.4mm,y=8mm]
  \node[case label] at (0,2) {$O_3$};
  \Comb{x/o,x/x,o/x}
  \LowerLabel{0}{v_1}
  \LowerLabel{1}{v_2}
  \UpperLabel{1}{w_1}
  \UpperLabel{2}{w_2}

  \begin{scope}[shift={(5.5,0)}]
    \node[case label] at (0,2) {$O_4$};
    \Comb{x/o,x/x,o/x,o/x}
    \LowerLabel{0}{v_1}
    \LowerLabel{1}{v_2}
    \UpperLabel{1}{w_1}
    \UpperLabel{2}{w_2}
    \UpperLabel{3}{w_3}
  \end{scope}

  \begin{scope}[shift={(0,-3.25)}]
    \CommonStart{0}
    \begin{scope}[shift={(9,0)}]
      \Comb{x/o,x/x,o/x,x/x,x/o,x/x,o/x,o/x}
      \LowerLabel{5}{v_{m-2}}
      \UpperLabel{6}{w_{m-1}}
      \UpperLabel{7}{w_m}
    \end{scope}
  \end{scope}

  \begin{scope}[shift={(0,-6.5)}]
    \CommonStart{1}
    \begin{scope}[shift={(9,0)}]
      \Comb{x/o,x/x,o/x,x/x,x/o}
      \LowerLabel{3}{v_{m-1}}
      \LowerLabel{4}{v_m}
      \UpperLabel{3}{w_{m-1}}
    \end{scope}
  \end{scope}

  \begin{scope}[shift={(0,-9.75)}]
    \CommonStart{2}
    \begin{scope}[shift={(9,0)}]
      \Comb{x/o,x/x,o/x,x/x,x/o,x/o}
      \LowerLabel{3}{v_{m-2}}
      \LowerLabel{4}{v_{m-1}}
      \LowerLabel{5}{v_m}
      \UpperLabel{3}{w_{m-2}}
    \end{scope}
  \end{scope}

  \begin{scope}[shift={(0,-13)}]
    \CommonStart{3}
    \begin{scope}[shift={(9,0)}]
      \Comb{x/o,x/x,o/x,x/x,x/o,x/x,o/x}
      \LowerLabel{5}{v_{m-1}}
      \UpperLabel{5}{w_{m-1}}
      \UpperLabel{6}{w_m}
    \end{scope}
  \end{scope}
\end{tikzpicture}}
\caption{Minimal super dominating sets of maximum cardinality of $O_{m}$ }\label{GammaspCm}
\label{Gammasp3}
\end{figure}

Now, we will prove that
$ \Gamma_{sp}(O_{m})
\leq \left\{
    \begin{array}{ll}
         3q+1, &if\, m=2q+1, q\geq 1,\\
        3q+2, &if\, m=2q+2, q\geq 1.
    \end{array}
    \right.$

Let $S_m$ be a minimal super dominating set of $O_m$.  Similarly like before, we consider some cases depending on the value of $m.$

\vspace{0.2cm}

{\bf Case1: $m=2q+1,$ with $q\geq 1.$}
Using Lemma \ref{lem:Cm} and Lemma \ref{lem:Cm2} we have that
\[
|S_{m}|=|S_{m}\cap \{v_1,w_1\}|+\sum_{i=1}^q|S_{m}\cap \{v_{2i},w_{2i},v_{2i+1},w_{2i+1}\}|\leq 1+3q.
\]
    
{\bf Case 2:  $m=2q+2,$ with $ q\geq 1$.} Using Lemma \ref{lem:Cm} and Lemma \ref{lem:Cm2} we have
\[
|S_{m}|=|S_{m}\cap \{v_1,w_1\}|+|S_{m}\cap \{v_m,w_m\}|+\sum_{i=1}^q|S_{m}\cap \{v_{2i},w_{2i},v_{2i+1},w_{2i+1}\}|
\]
$
\phantom{cccccll}\leq 1+1+3q=3q+2.
$

\vspace{0.2cm}

We conclude 
$$ \Gamma_{sp}(O_{m})
= \left\{
    \begin{array}{ll}
         3q+1, &if\, m=2q+1, q\geq 1,\\
        3q+2, &if\, m=2q+2, q\geq 1.
    \end{array}
    \right.$$
\end{proof}

\section{Some Gallai-type theorems}
In this section we present some Gallai-type theorems involving dominating, super dominating and open irredundant sets. 

\subsection{Relations between $\Gamma_{sp}(G)$ and $\gamma(G)$}
The calculations we showed in the previous section lead us to the following two propositions.

\begin{prop}\label{prop:Pm4}
Let $n\in\mathbb{N}^+$. Then, $\gamma(P_n) + \Gamma_{sp}(P_n) \geq n$.
\end{prop}

\begin{proof}
Let $q,r\in \mathbb{N}$ such that $n=21q+r$, with $0\leq r<21$. For $2\leq n\leq 4$, we have the equality $\gamma(P_n) + \Gamma_{sp}(P_n)=n$. For $n\geq 5$, we know that $\gamma(P_n)= \lceil\frac{n}{3}\rceil$ and by Proposition \ref{prop:Pm2} $\Gamma_{sp}(P_n)\geq 15q+\lceil\frac{2r-1}{3}\rceil$ so, 

$$\begin{array}{ll}\gamma(P_n) + \Gamma_{sp}(P_n) &\geq \lceil\frac{n}{3}\rceil+15q+\lceil\frac{2r-1}{3}\rceil=\lceil\frac{21q+r}{3}\rceil+15q+\lceil\frac{2r-1}{3}\rceil\\ &=7q+\lceil\frac{r}{3}\rceil+15q+\lceil\frac{2r-1}{3}\rceil=22q+\lceil\frac{r}{3}\rceil+\lceil\frac{2r-1}{3}\rceil.\end{array}$$

We consider some cases, depending on the value of $r$.

\vspace{0.2cm}

{\bf Case 1: $r=3\rho$, $\rho\in\mathbb{N}.$}
$$\begin{array}{ll}\gamma(P_n) + \Gamma_{sp}(P_n) &\geq 22q+\lceil\frac{r}{3}\rceil+\lceil\frac{2r-1}{3}\rceil=22q+\lceil\frac{3\rho}{3}\rceil+\lceil\frac{2(3\rho)-1}{3}\rceil\\&=22q+\rho+2\rho=q+(21q+3\rho)=q+n\geq n.\end{array}$$

{\bf Case 2: $r=3\rho+1$, $\rho\in\mathbb{N}.$}
$$\begin{array}{ll}\gamma(P_n) + \Gamma_{sp}(P_n) &\geq 22q+\lceil\frac{r}{3}\rceil+\lceil\frac{2r-1}{3}\rceil=22q+\lceil\frac{3\rho+1}{3}\rceil+\lceil\frac{2(3\rho+1)-1}{3}\rceil\\&=22q+\rho+1+2\rho+1=q+(21q+3\rho+1)+1=q+n+1> n.\end{array}$$

{\bf Case 3: $r=3\rho+2$, $\rho\in\mathbb{N}.$}
$$\begin{array}{ll}\gamma(P_n) + \Gamma_{sp}(P_n) &\geq 22q+\lceil\frac{r}{3}\rceil+\lceil\frac{2r-1}{3}\rceil=22q+\lceil\frac{3\rho+2}{3}\rceil+\lceil\frac{2(3\rho+2)-1}{3}\rceil\\&=22q+\rho+1+2\rho+1=q+(21q+3\rho+2)=q+n\geq n.\end{array}$$
We conclude that $\gamma(P_n) + \Gamma_{sp}(P_n) \geq n$.

\end{proof}

\begin{obs}\label{obs:Pn4} Notice that there exist $n\in\mathbb{N}^+$ such that $\gamma(P_n) + \Gamma_{sp}(P_n) > n$, for example $n=21q+r$, with $q,r\in \mathbb{N}$, $0\leq r<21$ and $q>0$. 
\end{obs}

The following proposition shows the strict inequality fulfill for caterpillars:

\begin{prop}\label{Cm2}
Let $m\in\mathbb{N}^+$ with $m>2$. Then, $\gamma(O_{m}) + \Gamma_{sp}(O_m) > 2m$.
\end{prop}

\begin{proof}
Let $m\in\mathbb{N}^+$ with $m>2$. Consider the following cases:

\vspace{0.2cm}

{\bf Case1:   $m=2q+1,$ with $q\geq 1.$}
By Proposition \ref{Cm}, $\Gamma_{sp}(O_{m})=3q+1$. Also, $D_{m}=\{v_j\,|\,j=1,\dots,m\}$ is a dominating set of $O_m$ of minimum cardinality, so $\gamma(O_{m})=|D_{m}|=m$. Then, $$\gamma(O_{m}) + \Gamma_{sp}(C_{m})=m+3q+1=2q+1+3q+1=6q+2>4q+2=2m.$$ 

\vspace{0.2cm}

{\bf Case 2: $m=2q+2,$ with $ q\geq 1$.}
By Proposition \ref{Cm}, $\Gamma_{sp}(O_{m})=3q+2$. Also, $D_{m}=\{v_j\,|\,j=1,\dots,m\}$ is a dominating set of $O_m$ of minimum cardinality, so $\gamma(O_{m})=|D_{m}|=m$. Then, $$\gamma(O_{m}) + \Gamma_{sp}(O_{m})\geq m+3q+2=2q+2+3q+2=6q+4>4q+4=2m.$$

\end{proof}

In 1979 Bollob\'{a}s and Cockayne \cite{Bollobas79} proved the following now well-known result.

\begin{thm} (Bollob\'{a}s-Cockayne \cite{Bollobas79})\label{Bollobas}
Every isolate-free graph $G$ contains an open irredundant $\gamma$-set, that is a minimum cardinality dominating set $S$ in which every vertex $u \in S$ has at least one external private neighbour. 
\end{thm}
  With this in mind we can prove the following result, which will be useful to establish a relation between the domination number $\gamma(G)$ and the upper super domination number $\Gamma_{sp}(G)$.

\begin{lemma}\label{ten}
If $D$ is a $\gamma$-set in a graph $G$ that is open irredundant, then $V - D$ is a minimal super dominating set of $G$.
\end{lemma}

{\bf Proof}.  Let $D$ be a $\gamma$-set in a graph $G$ that is open irredundant. Thus, $D$ is a private dominating set and from Proposition~\ref{six}, we have that $V - D$ is a super dominating set. Assume that $V - D$ is not a minimal super dominating set; thus there exists a vertex $x \in V - D$ such that $(V - D) - \{x\}$ is a super dominating set. 

Hence, from Proposition~\ref{six}, $D \cup \{x\}$ is an open irredundant set.  This means that vertex $x$ has an external private neighbour with respect to $D \cup \{x\}$.  Thus, there exists a vertex $z \in V - (D \cup \{x\}) = (V - D) - \{x\}$ which is an external private neighbour of $x$ with respect to $D \cup \{x\}$.  In particular, vertex $z$ has no neighbours in $D$, which contradicts the fact that $D$ is is a $\gamma$-set of $G$. Thus, $V - D$ is a minimal super dominating set of $G$. $\Box$

Now we can generalize Proposition \ref{prop:Pm4}.
\begin{prop}\label{eleven} 
For any graph $G$ of order $n$, $\gamma(G) + \Gamma_{sp}(G) \ge n$.
\end{prop}

\begin{proof}  Let $D$ be a $\gamma$-set in a graph $G$ such that $D$ is open irredundant (from Theorem~\ref{Bollobas} we know that such a $\gamma$-set exists). By Lemma~\ref{ten}, $V - D$ is a minimal super dominating set, so $|V - D| \le \Gamma_{sp}(G)$. Thus, we have that $n - \gamma(G) \le \Gamma_{sp}(G)$, or equivalently, $\Gamma_{sp}(G) + \gamma(G) \ge~n$. \end{proof}

Now we give some conditions that guarantee the equality in the inequality of Proposition \ref{eleven}.

\begin{prop}
If there exists a $\Gamma_{sp}$-set of $G$ that is also a minimal dominating set of $G$, then $\gamma(G) + \Gamma_{sp}(G) = n$.
\end{prop}

\begin{proof}  Let $G$ be a graph having a  $\Gamma_{sp}$-set $D$ that is also a minimal dominating set of $G$.  From Proposition~\ref{eleven} we have that $\gamma(G) + \Gamma_{sp}(G) \ge n$.  From Observation~\ref{seven}, $V - D$ is an open irredundant set. In \cite{Ore62} Ore proved that the complement of a minimal dominating set is a dominating set, so $V - D$  is also a dominating set.  Thus, $|V - D| \ge \gamma(G)$, so $n - \Gamma_{sp}(G) \ge \gamma(G)$.  Thus, $\gamma(G) + \Gamma_{sp}(G) \le n$, and finally, $\gamma(G) + \Gamma_{sp}(G) = n$. \end{proof}

\begin{cor} \label{thirteen}
If there exists a $\Gamma_{sp}$-set $D$ in a graph $G$ that is a dominating set and its complement $V - D$ is also a dominating set, then $\gamma(G) + \Gamma_{sp}(G) = n$.
\end{cor}

\begin{obs}
If a graph $G$ has a universal vertex, then $\gamma(G) = 1$,  $\Gamma_{sp}(G) = n - 1$, and $\gamma(G) + \Gamma_{sp}(G) = n$.
\end{obs}

In \cite{Detlaff19}, the authors give a characterization of  graphs $G$ such that $\gamma _{sp}(G)=n-1, \gamma (G)=1$,  which satisfy the observation below.

\subsection{Super domination and enclaveless parameters}

In 1977 Slater \cite{Slater77} defined an {\it enclave} in a set $S \subseteq V$ to be a vertex $v \in S$ such that $N[v] \subseteq S$, that is, $v$ and all of its neighbours are in $S$.  A set S is {\it enclaveless} if it contains no enclaves.  One can observe that a set $S$ is a dominating set of a graph G if and only if its complement $V - S$ is enclaveless.  The maximum cardinality of an enclaveless set in a graph G is denoted by $\Psi(G)$ and is called the {\it enclaveless number} of $G$, while the {\it lower enclaveless number} $\psi(G)$ equals the minimum cardinality of a maximal enclaveless set. In \cite{Slater77} Slater established the following result.

\begin{thm} (\cite{Slater77})\label{Slater}
For any graph $G$ of order $n$, $\gamma(G) + \Psi(G) = n$.
\end{thm}

The simple proof of this theorem immediately yields the following two corollaries.

\begin{cor}
For any graph $G$, the complement of any (minimal) dominating set is a (maximal) enclaveless set, and the complement of any (maximal) enclaveless set is a (minimal) dominating set.
\end{cor}

\begin{cor}
For any graph $G$ of order $n$, $\Gamma(G) + \psi(G) = n$.
\end{cor}

Thus, from Proposition~\ref{eleven} and Theorem~\ref{Slater}, we have the following inequality.

\begin{prop}
For any graph $G$, $\Gamma_{sp}(G) \ge \Psi(G)$.
\end{prop}

\subsection{Super domination and open irredundance parameters}

We begin with relations between open irredundant and super dominating sets.

\begin{prop}\label{six}
$S \subseteq V$ is an open irredundant set of $G$ if and only if the set $V - S$ is a super dominating set of $G$.
\end{prop}

\begin{proof}  A set $S$ is an open irredundant set of $G$ if and only if every vertex $u \in S$ has an external private neighbour with respect to $S$.  This implies that $V - S$ dominates the set $S$ and for every vertex $u \in V - S$ there is a vertex $v \in V - (V - S) = S$ whose only neighbor in $V - S$ is $u$, and thus $V - S$ is a super dominating set of $G$. \end{proof}

In addition, the following is true.

\begin{obs}\label{seven}
If $S \subseteq V$ is a super dominating set of $G$, then $V - S$ is an open irredundant set, but not necessarily a dominating set. 
\end{obs}

For example, any four consecutive vertices of a cycle $C_6$ form a (minimal) super dominating set whose complement fails to dominate $C_6.$

As a consequence of Proposition~\ref{six} we can relate minimal super dominating sets of $G$ with maximal open irredundant sets of $G$.

\begin{prop}\label{propeight}
A set $S \subseteq V$ is a maximal open irredundant set of $G$ if and only if $V - S$ is a minimal super dominating set of $G$.
\end{prop}

\begin{proof}  Let $S$ be a maximal open irredundant set.   From Proposition~\ref{six} it follows that $V - S$ is a super dominating set of $G$.  Assume that $V - S$ is not a minimal super dominating set. Then there exists a super dominating set $T \subset V - S$.  By Proposition~\ref{six} we know that $V - T$ is an open irredundant set, but $S \subset V - T$, which contradicts the assumption that $S$ is a maximal open irredundant set of $G$. Thus, $V - S$ is a minimal super dominating set.

Conversely, let $D$ be a minimal super dominating set of $G$.  From Proposition~\ref{six}, $V - D$ is an open irredundant set.  Assume that $V - D$ is not a maximal open irredundant set.   Then there exists an open irredundant set $F$ such that $V - D \subset F$. By Proposition~\ref{six} we know that $V - F$ is a super dominating set of $G$. But $V - F \subset D$, which contradicts the fact that $D$ is a minimal super dominating set of $G$. Thus, $V - D$ is a maximal open irredundant set. \end{proof}

Recall from the Introduction that $oir(G)$ and $OIR(G)$ equal the minimum and maximum cardinality of a maximal open irredundant set, respectively. In this section we establish two Gallai theorems involving super domination and open irredundance.

\begin{thm}\label{twentytwo}
For any graph $G$ of order $n$, $\gamma_{sp}(G) + OIR(G) = n$.
\end{thm}

\begin{proof}  Let $D$ be a minimum cardinality super dominating set of a graph $G$, so that $|D| = \gamma_{sp}(G)$.  Then $|D| \le |S|$ for every super dominating set $S$ of $G$.  So, $-|D| \ge -|S|$ for every super dominating set $S$ of G and therefore, $|V - D| = n - |D| \ge n - |S| = |V - S|$ for every super dominating set $S$ of $G$.  By 
Proposition~\ref{six}, this is equivalent to the fact that $|V - D| = n - |D| \ge |T|$ for every open irredundant set $T$ of $G$. Then $V - D$ is an open irredundant set of $G$ of maximal cardinality, so $|V - D| = OIR(G)$.  Therefore, $n = |D| + |V - D| = \gamma_{sp}(G) + OIR(G)$. \end{proof}

\begin{thm}\label{twentythree}
For any graph $G$ of order $n$, $\Gamma_{sp}(G) + oir(G) = n$.
\end{thm}

\begin{proof} Let $D$ be a minimal super dominating set in a graph $G$ of maximum cardinality, so that $D = \Gamma_{sp}(G)$. Then $|D| \ge |S|$ for every minimal super dominating set $S$ of $G$, so $-|D| \le -|S|$, for every minimal super dominating set $S$ of $G$.  Therefore,  $|V - D| = n - |D| \le n - |S| = |V - S|$ for every minimal super dominating set $S$ of $G$.  By Proposition~\ref{propeight} this implies that $|V - D| = n - |D| \le |T|$ for every maximal open irredundant set $T$ of $G$. Therefore, $V - D$ is a maximal open irredundant set of $G$ of minimum cardinality, and so, $|V - D| = oir(G)$. Thus, $n = |D| + |V - D| = \Gamma_{sp}(G) + oir(G)$. \end{proof}

From Theorem~\ref{twentytwo}, Theorem~\ref{twentythree}, and Proposition~\ref{propGammasp}, we can conclude that the difference $OIR(G) - oir(G)$ can be arbitrarily large.

\begin{cor}
For any positive integer $k \in \mathcal{N}$ there exists a graph $G$ such that $OIR(G) - oir(G) = k$.
\end{cor}

\section{Aknowledgements}
We would like to express our sincere thanks to Professor Stephen Hedetniemi for his kind support and for the inspiration to write this article.

\end{document}